\documentclass[letterpaper,oneside,10pt]{article}

\usepackage[USenglish]{babel} 
\usepackage[T1]{fontenc}
\usepackage[ansinew]{inputenc}

\usepackage{lmodern} 

\usepackage{graphicx} 

\usepackage{amsmath}
\usepackage{amsthm}

\newtheorem{lemma}{Lemma}
\usepackage{amsfonts}
\usepackage{amssymb}
\usepackage{mathtools}
\usepackage[hyphens]{url}

\let\Right\right 
\let\Left\left 
\makeatletter 
\def\right#1{\Right#1\@ifnextchar){\!\right}{}} 
\def\left#1{\Left#1\@ifnextchar({\!\left}{}} 
\makeatother

\begin{document}

\pagestyle{empty} 


\title{On the Andrews congruence for the Fibonacci quotient}
\author{John Blythe Dobson (j.dobson@uwinnipeg.ca)}

\maketitle


\pagestyle{plain} 


\begin{abstract}
\noindent
We show that a congruence discovered by George E.\ Andrews in 1969 for the Fibonacci quotient directly implies a simpler congruence found by Hugh C.\ Williams in 1991.

\noindent
\textit{Keywords}: Fibonacci quotient, harmonic numbers
\end{abstract}

\section{Introduction}

\noindent
In his celebrated 1969 study of the Fibonacci sequence (\cite{Andrews}, p.\ 114), George E.\ Andrews gave two congruences modulo a prime $p$ for divided Fibonacci numbers:

\begin{subequations}

\begin{equation} \label{eq:Andrews}
F_{p-1}/p \equiv 2 \left(\frac{-1}{\,p}\right) \sum_{\substack{j=-p+1\\ j \equiv 5, 7 \bmod{10}}}^{p-1} \frac{\left(\frac{j+1}{5}\right)\left(\frac{-1}{j}\right)}{p-j} \quad [p \equiv \pm 1 \pmod{5}; p > 5]
\end{equation}

\begin{equation}
F_{p+1}/p \equiv 2 \left(\frac{-1}{\,p}\right) \sum_{\substack{j=-p+1\\ j \equiv 1, 5 \bmod{10}}}^{p-1} \frac{\left(\frac{j+1}{5}\right)\left(\frac{-1}{j}\right)}{p-j} \quad [p \equiv \pm 2 \pmod{5}; p > 5].
\end{equation}

\end{subequations}

\noindent
Andrews takes the Jacobi symbol $\left(\frac{-1}{\,p}\right)$ in the sense $(-1)^{(p - 1)/2}$, while some authors take it as $(-1)^{(\lvert p \rvert - 1)/2}$. Together these congruences characterize the Fibonacci quotient $F_{p-(\frac{5}{p})}/p$ for $p \neq 5$, since $\left(\frac{5}{p}\right)$ = 1 and $-1$ in the two cases, respectively. Andrews noted the resemblance of the sums in these results with the sums of reciprocals that occur in the congruence for the Fermat quotient in the classic study of Eisenstein \cite{Eisenstein}. Later Hugh C.\ Williams, using (as he himself notes) a quite different method, was able in his 1982 paper on the Fibonacci quotient \cite{Williams1982} to derive a simpler congruence exactly analogous to Eisenstein's, and in his 1991 paper treating quotients of Lucas numbers (\cite{Williams1991}, p.\ 440, eq.\ 4.7) he further simplified this to

\begin{equation} \label{eq:Williams}
F_{p-(\frac{5}{p})}/p \equiv \frac{2}{5} \sum_{j = \lfloor p/5 \rfloor + 1}^{\lfloor 2p/5 \rfloor} \frac{1}{j} \pmod{p} \quad [p > 5].
\end{equation}

\noindent
Willliams's result can also be viewed as a statement about a special type of harmonic number $H$, since it can be written

\begin{equation} \label{eq:WilliamsHarmonic}
F_{p-(\frac{5}{p})}/p \equiv \frac{2}{5} \left\{ H_{\lfloor 2p/5 \rfloor} - H_{\lfloor p/5 \rfloor} \right\} \pmod{p} \quad [p > 5].
\end{equation}

\noindent
We regard Williams's elegant formulation as canonical, as it is the simplest possible when the terms in the sum are constrained to be of like sign and equal weight. This property facilitates comparison with other results pertaining to the Fermat quotient, particularly those discussed in \cite{Lehmer1938} and \cite{DilcherSkula1995} in relation to the first case of Fermat's Last Theorem.

As it does not appear to have been generally recognized that Andrews's result implies (\ref{eq:Williams}), and as the exercise of deriving (\ref{eq:Williams}) directly from (\ref{eq:Andrews} \& b) may perhaps suggest possibilities for the simplification of similar results in the literature, we shall demonstrate how this can be done. We do so with unavoidable foreknowledge of (\ref{eq:Williams}), but without introducing any additional facts about the Fibonacci numbers. The only apparatus required follows from a 1905 paper of Lerch \cite{Lerch}.

\section{Some preliminaries}

\noindent
We shall make use of a notation originally introduced by Lerch, but with the shift in index adopted in \cite{DilcherSkula1995} and most recent writings:

\begin{equation} \label{eq:SkulaDefinition}
s(k, N) = \sum_{\substack{j=\lfloor\frac{kp}{N}\rfloor + 1\\ j \neq p}}^{\lfloor\frac{(k + 1)p}{N}\rfloor}
\frac{1}{j},
\end{equation}

\noindent
where it is always assumed that $p$ is sufficiently large that $s(k, N)$ contains at least one element; the provision $j \neq p$ is necessary when $k+1=N$, though we shall not encounter that situation here. In this notation, therefore, the sum in the right-hand side of (\ref{eq:Williams}) can be written as $s(1, 5)$. We shall frequently make use of the trivial fact that $s(k, N) \equiv -s(N-1-k, N)$ $\pmod{p}$.

Lerch (\cite{Lerch}, p.\ 476, equations 14 and 15), correcting work of Sylvester, pretty much explicitly writes out the relation $2 \cdot s(0, 5) + s(1, 5) \equiv -\frac{5}{2} \cdot q_p(5)$, where $q_p(5)$ is the Fermat quotient of $p$ to the base 5. While we shall not pursue this matter here, it is thus evident that the evaluation of $s(1, 5)$ simultanously settled the evaluation of $s(0, 5)$. This result, and the core formulae of Lerch's paper, supply a number of relationships between sums of the type defined in (\ref{eq:SkulaDefinition}), and these are greatly extended in a recent paper of Dilcher and Skula \cite{DilcherSkula2011}.

We shall need a lemma in this vein:

\begin{lemma} \label{Lemma 1}
\begin{equation} \label{eq:Lemma1}
-s(1, 10) + 2 \cdot s(2, 10) + 3 \cdot s(3, 10) \equiv 0 \pmod{p}.
\end{equation}

\begin{proof}
The ingredients needed for this lemma are developed using only ideas from Lerch \cite{Lerch} in \cite{DilcherSkula2011}, pp. 20--21, Proposition 3.1. Earlier, it had been proved in a less elementary way, and employing a slightly different notation, in a paper of Skula (\cite{Skula2008}, p.\ 9, Theorem 3.2), which gives:

\begin{displaymath}
2 \cdot s(0, 10) + 3 \cdot s(1, 10) + 2 \cdot s(2, 10) + 3 \cdot s(3, 10) + 2 \cdot s(4, 10) \equiv 0 \pmod{p}.
\end{displaymath}

\begin{displaymath}
s(0, 10) + 2 \cdot s(1, 10) + s(4, 10) \equiv 0 \pmod{p}.
\end{displaymath}

\noindent
Subtracting twice the second row from the first gives (\ref{eq:Lemma1}). 
\end{proof}

\end{lemma}

We shall also need a formula for a sum which, under various notations, has made frequent appearances in the literature of the Fermat quotient. If $K(r, N)$ represents the sum of the terms in $s(0, 1)$, i.e. the sum of the terms in the set $\left\{ 1, \frac{1}{2}, \frac{1}{3}, \dots, \frac{1}{p-1} \right\}$, whose denominators are congruent to $rp \pmod{N}$ for a prime $p$ and $r < N$, then

\begin{lemma} \label{Lemma 2}

\begin{equation} \label{eq:Lemma2}
K(r, N) \equiv \frac{1}{N} \cdot s(N-r, N) \pmod{p}.
\end{equation}

\begin{proof}
This formula was given in 1995 by Zhi-Hong Sun (\cite{ZHSun1992}, pt. 3, p. 90, Corollary 3.1), though it may be older. For an elementary proof using ideas from Lerch \cite{Lerch} see \cite{DobsonLerchFormula}, \S{}4.
\end{proof}

\end{lemma}

\noindent
This relation, which defines an association between terms characterized by a congruential condition on denominators lying in the interval $\left\{ 1, p-1 \right\}$, and a consecutive block of terms $s(k, N)$, permits simplification of many published results involving sums of reciprocals. One is usually interested in sums of terms whose denominators belong to some fixed residue class $t$, so if $p$ is invertible modulo $N$ and this inverse is $p^\prime$, one sets $r = tp^\prime$.

\section{The main result}

\noindent
We are now ready to derive (\ref{eq:Williams}) from (\ref{eq:Andrews} \& b). First, we rewrite (\ref{eq:Andrews} \& b) in a more explicit form, still modulo $p$ but with the summations now confined to the range $\left\{ 1, p-1 \right\}$:

\begin{subequations}

\begin{equation} \label{eq:AndrewsRewritten}
\begin{split}
& F_{p-1}/p \\
& \equiv 2 \left(\frac{-1}{\,p}\right) \left\{ \sum_{\substack{j \equiv 15\\ \bmod{20}}} \frac{2}{j} \quad - \sum_{\substack{j \equiv 5 \\ \bmod{20}}} \frac{2}{j} \quad + \sum_{\substack{j \equiv 13, 17 \\ \bmod{20}}} \frac{1}{j} \quad - \sum_{\substack{j \equiv 3, 7 \\ \bmod{20}}} \frac{1}{j} \right\} \\
& \qquad \qquad \qquad \qquad \qquad \qquad \qquad \qquad \qquad \qquad \qquad [p \equiv \pm 1 \bmod{5}; p > 5]
\end{split}
\end{equation}

\begin{equation}
\begin{split}
& F_{p+1}/p \\
& \equiv 2 \left(\frac{-1}{\,p}\right) \left\{ \sum_{\substack{j \equiv 15\\ \bmod{20}}} \frac{2}{j} \quad - \sum_{\substack{j \equiv 5 \\ \bmod{20}}} \frac{2}{j} \quad + \sum_{\substack{j \equiv 1, 9 \\ \bmod{20}}} \frac{1}{j} \quad - \sum_{\substack{j \equiv 11, 19 \\ \bmod{20}}} \frac{1}{j} \right\} \\
& \qquad \qquad \qquad \qquad \qquad \qquad \qquad \qquad \qquad \qquad \qquad [p \equiv \pm 2 \pmod{5}; p > 5]
\end{split}
\end{equation}

\end{subequations}

\noindent
The first row covers the cases $p \equiv 1, 9, 11, 19 \pmod{20}$, and the second row the cases $p \equiv 3, 7, 13, 17 \pmod{20}$; the eight cases must be split out in order to determine the values of the Jacobi symbol $\left(\frac{-1}{\,p}\right)$. Next, we apply Lemma \ref{Lemma 2} to each of the component sums, letting $t$ be the residue class modulo 20 specified for each summation, and for $p \equiv$ 1, 9, 11, 19, 3, 7, 13, 17 $\pmod{20}$, taking $p^\prime$ = 1, 9, 11, 19, 7, 3, 17, 13, respectively (the first four residues are their own inverses, being the square roots of unity). Routine calculations then establish that for all eight cases of $p$, (\ref{eq:AndrewsRewritten} \& b) reduce (after rearrangement) to

\begin{equation} \label{eq:AndrewsSimplified}
\begin{split}
& F_{p-(\frac{5}{p})}/p \\
& \equiv \frac{1}{10} \left\{ s(2, 20) + s(3, 20) + 2 \cdot s(4, 20) + 2 \cdot s(5, 20) + s(6, 20) + s(7, 20) \right\} \\
& \qquad \qquad \qquad \qquad \qquad \qquad \qquad \qquad \qquad \qquad \qquad \pmod{p} \quad [p > 5].
\end{split}
\end{equation}

\noindent
In view of the definition of $s(k, N)$, it is clear that when $k$ and $N$ are both even, we have $s(k, N)$ + $s(k+1, N)$ = $s(\frac{k}{2}, \frac{N}{2})$. Thus (\ref{eq:AndrewsSimplified}) condenses to

\begin{displaymath}
F_{p-(\frac{5}{p})}/p \equiv \frac{1}{10} \left\{ s(1, 10) + 2 \cdot s(2, 10) + s(3, 10) \right\} \pmod{p} \quad [p > 5],
\end{displaymath}

\noindent
and finally, adding one tenth of (\ref{eq:Lemma1}) to this gives

\begin{displaymath}
F_{p-(\frac{5}{p})}/p \equiv \frac{1}{10} \left\{ 4 \cdot s(2, 10) + 4 \cdot s(3, 10) \right\} \equiv \frac{2}{5} \cdot s(1, 5) \pmod{p} \quad [p > 5],
\end{displaymath}

\noindent
which is equivalent to (\ref{eq:Williams}), as required.

\end{document}